\documentclass[11pt,a4paper]{article}

\usepackage{latexsym,amsfonts,amstext,amsbsy,amsmath,amssymb,amscd,graphics,epsfig,eucal,bbm}

\usepackage{pb-diagram}
\usepackage[all]{xy}
\usepackage{makeidx}
\usepackage{indentfirst}
\usepackage{graphicx}
\usepackage{mathrsfs}

\begin{document}
\begin{center}
{\large \bf {A TESTING ALGORITHM OF AN UNIVERSAL ALGEBRA TO BE A
BRANDT GROUPOID }}
\end{center}
\begin{center}
{\bf Gheorghe IVAN and George STOIANOV}
\end{center}
\date{}
\def \e{\varepsilon}
\def \w{\widetilde}
\def\b{\Box}
\begin{center}
{\bf Abstract.}
\end{center}
The main aim of this paper is to present a program on computer for
decide if an universal algebra is a groupoid. Using the theory of
groupoids and the program $~BGroidAP1~$ we prove a theorem of
classification for the groupoids of type $~(4;2).~$
{
\begin{center}
{\bf  Introduction}
\end{center}

The algebraic notion of groupoid was introduced and named by H.
Brandt in the paper [ {\it Uber eine Verallgemeinerung der
Gruppen-begriffes}. Math. Ann., 96, 1926, 360-366 ]. A groupoid (
in the sense of Brandt ) can be thought of as a generalized group
in which only certain multiplications are possible and it contains
several neutral elements.

Groupoids also appeared in Galois theory in the description of
relations between subfields of a field $~K~$ via morphisms of
$~K~$ in a paper of A. Loewy [ {\it Neue elementare Begrundung und
Erweiterung der Galoisschen Theorie}. S.-B. Heidelberger Akad.
Wiss. Math. Nat. Kl.{\bf 1925}, 1927 ]. In differentiable context,
the concept of groupoid has appeared  in the work of C. Ehresmann
[{\it Cat$\acute e$gories et structures}. Dunod, Paris ] around
1950.

In the language of categories, a groupoid is a small category in
which all morphisms are invertible. For more details concerning
the groupoids and its applications in many areas of mathematics,
see [1] - [4],[6] - [8].

The plan of this paper is as follows. In the first section we have
collected the preliminary concepts concerning groupoids. In the
second section we present an algorithm for to verify that a finite
set endowed with structure functions has a groupoid structure.
This algorithm is based on the theory of groupoids and is
implemented on computer. The obtained program is denoted by
$~BGroidAP1.~$ Finally, we ilustrate the utilisation of the
program on some finite universal algebras. In particular, the
program can be used for to test if a finite set endowed with a
composition law has a structure of group.

\begin{center}
{\bf 1. The concept of Brandt groupoid as universal algebra}
\end{center}

{\bf Definition 1.1.}~Let $~( M,M_{0})~$ be a pair of nonempty
sets, where $~M_{0}~$ a subset of $~M~$ endowed with the
surjections $~\alpha, ~\beta : M~\to~M_{0}~$, called  the {\it
source}  and the {\it target} map, respectively and a {\it (
partial ) multiplication  law} $~\mu :M_{(2)}\longrightarrow M,
(x,y)~\longrightarrow ~\mu(x,y),~$ where $~M_{(2)} =\{~(x,y)\in M
\times M~|~\beta (x)=\alpha (y)~\}.~$ We write sometimes $~x\cdot
y~$ or $~x y~$ for $~\mu(x,y).~$ The elements of $~M_{(2)}~$ are
called {\it composable pairs}~ of $~M.~$

We say that the universal algebra $~( M, \alpha, \beta, \mu ;
M_{0})~$ is a {\it semigroupoid}, if the multiplication law is
{\it associative}, i.e. $~ (x y) z = x (y z),~$ in the sense that,
if one side of the equation is defined so is the other one and
then they are equal( the element $~(x y) z~$ is defined iff
$~\beta(x) = \alpha(y)~$ and $~\beta(y) = \alpha(z)).~$\hfill$\b$

{\bf Definition 1.2.}~([4]) A {\it monoidoid} is a semigroupoid
$~(M, \alpha, \beta, \mu;M_{0} )~$ such that the {\it identities
property} holds, i. e.  for each $~x\in M~$ we have $~(\alpha
(x),x),~(x,\beta (x))\in M_{(2)}~$ and $~ \alpha(x) x = x\beta(x)
= x.$ \hfill$\b$

The element $~\alpha(x)~$ [ resp. $~\beta(x)~$ ] denoted sometimes
by $~u_{l}(x)~$ [ resp. $~u_{r}(x)~$ ] is the {\it left unit}  [
resp. {\it right unit}]  of $~x\in M.~$ The subset $~M_{0} =
\alpha(M) = \beta(M)~$ of $~M~$  is called the {\it unit set} of
$~M~$ and we say that $~M~$ is a $~M_{0}~$ - {\it monoidoid}. The
functions $~\alpha, \beta, \mu~$ are called {\it structure
functions} of the monoidoid $~M.~$ For each $~u\in M_{0},~$ the
set $~\alpha^{-1}(u)~$  ( resp., $~\beta^{-1}(u)~$ ) is called
$~\alpha~$-{\it fibre} ( resp., $~\beta~$-{\it fibre} ) of the
monoidoid $~M~$ over $~u\in M_{0}.~$

In the following proposition we summarize some properties of the
structure functions of a monoidoid.

{\bf Proposition 1.1.}~( [4] )~{\it Let $~(M, \alpha, \beta, \mu; M_{0})~$ and $~u,v\in M_{0}.~$ Then the following assertions hold:\\
$(1.1)~~~\alpha(u) = \beta(u) = u ~$ and $~u\cdot u = u~$ for all $~u\in M_{0};~$\\
$(1.2)~~~\alpha(x y) = \alpha(x)~$ and $~ \beta(x y ) = \beta(y), ~(\forall)~ (x,y)\in M_{(2)};~$\\
$(1.3)~~~$ if  $~(x,u)\in M_{(2)}~$ such that $~x u = x~$ then $~u=\beta(x);~$\\
$(1.4)~~~$ if  $~(v,x)\in M_{(2)}~$ such that $~v x = x~$ then
$~v=\alpha(x);~$}\hfill$\b$

{\bf Example 1.1.}~(i)~ A monoid $~{\cal M}~$ having $~e~$ as
unity, is just a $~\{ e\}~$- monoidoid over a one point-base in
the following way: $~M={\cal M}, M_{0} = \{e\},~$ the source and
target maps $~\alpha, \beta : M~\to~M_{0}~$ are constant maps,
i.e. $~\alpha(x)=\beta(x)=e~$ for all $~x\in M~$; for all $~x,y\in
M~$ we have $~\alpha(x)=\beta(x)=e~$ and hence the product
$~x\cdot y~$ is always defined in $~M~$ ( $~x\cdot y~$ is the
product of elements $~x~$ and $~y~$ in the monoid $~{\cal M}~$).

Conversely, every monoidoid $~M~$ with one unit ( i.e. $~M_{0}~$
is a singleton ) is a monoid.

$(ii)~$ The {\it nul monoidoid over a set}. Any nonempty set $~X~$
may be regarded as a monoidoid on itself with the following
monoidoid structure : $~M = M_{0} = X, \alpha = \beta = Id_{X}~$ ;
the elements $~x , y\in X~$ are composable iff $~ x=y~ $ and we
define $~ x\cdot x = x.~ $

$(iii)~$ {\it Disjoint union of two monoidoids}. Let
$~(M_{i},\alpha_{i}, \beta_{i}, \mu_{i}; M_{0,i})~, i = 1,2~$ be
two monoidoids such that $~M_{1}\cap M_{2} = \emptyset.~$ We
consider $ ~M = M_{1}\cup M_{2}~$ and $~M_{(2)}= M_{1,(2)}\cup
M_{2,(2)}~$.

We give on the set $~M~$ the following structure:\\
$~\alpha(x)=\alpha_{1}(x)~$ if $~x\in M_{1}~$ and $~\alpha(x) = \alpha_{2}(x)~$ if $~x\in M_{2};~$\\
$~\beta(x)=\beta_{1}(x)~$ if $~x\in M_{1}~$ and $~\beta(x) = \beta_{2}(x)~$ if $~x\in M_{2};~$\\
we have that $~(x,y)~\in M_{(2)}~$ iff $~(x,y)\in M_{1,(2)}~$ or
$~(x,y)\in M_{2,(2)}~$ and we take $~\mu(x,y) = \mu_{i}(x,y)~$ if
$~(x,y)\in M_{i,(2)}~, i=1,2.$

In other words, two elements $~x,y \in M$ may be composed iff they
lie in the same monoidoid $~M_{i}, i = 1,2~$ and they are
composable in $~ M_{i}, i = 1,2$.

In the case when $~M_{1}\cap M_{2}\neq \emptyset,~$ we cosider the
sets $~M_{1}^{\prime}=M_{1}\times \{1\},~ M_{2}^{\prime} =
M_{2}\times \{2\}~$ and we give on the set $~M_{1}^{\prime}\cup
M_{2}^{\prime}~$ the above monoidoid structure.

This monoidoid is denoted by $~M_{1}\coprod M_{2}~$ and is called
the {\it disjoint union of monoidoids} $~M_{1}~$ and $~M_{2}.~$
Its unit set is $~ M_{0} = M_{1,0}\cup M_{2,0}$, where $ M_{i,0}$
is the unit set of $ M_{i}, ~i=1,2$.

In particular, the disjoint union of monoids $~M_{i},~ i=1,2~$ is
a monoidoid, called the {\it monoidoid associated to monoids}
$~M_{i},~ i=1,2~$ ( for this monoidoid, the unit set
 is $~M_{0}= \{ e_{1}, e_{2}\}~$ where $~e_{i}~$ is the unity of $~M_{i}, i = 1,2~$ ).

$(iv)~$ Let the nul monoid $~M_{1} =\{e\}~$ and the multiplicative
monoid $~M_{2}=\{~-1, 0, 1~\}\subset {\bf Z}.~$ Then $~M=
M_{1}\cup M_{2}~= ~\{~e, 1, 0, -1~\}~$ is a monoidoid over $~M_{0}
= \{e,1~\}.~$\hfill$\b$

{\bf Example 1.2.}~$~(i)~$ {\it The monoidoid $~{\cal F}(S,X)~$}.
Let $~X~$ be a given nonepmty set. We denote by ${\cal F}(S,X) =
~\{~f : S~\to~X~|~ (\forall)~S~\hbox{such that}~\emptyset\neq
S\subseteq X~\}.~$ For $~f\in {\cal F}(S,X),~$ let $~D(f)~$ be the
domain of $~f~$ and let $~R(f) = f(D(f)).~$  For $~M={\cal
F}(S,X)~$ let $~M_{(2)}=\{(f,g)\in M\times M~|~R(f)=D(g)~\}~$ and
for $~(f,g)\in M_{(2)}~$ define $~\mu(f,g)=g\circ f.~$ If
$~Id_{S}~$ denotes the identity map on $~S,~$ then $~M_{0} =
\{~Id_{S}~|~\emptyset \neq S\subseteq X~\}~$ is the set of units
of $~M.~$ The maps $~\alpha, \beta : M~\to~M~$ are defined by
$~\alpha(f) = Id_{D(f)},~\beta(f)=Id_{R(f)}.~$  Thus $~{\cal
F}(S,X)~$ ia a monoidoid, called the {\it monoidoid of functions
from $~S~$ to $~X~$}, where $~S~$ is an arbitrary nonempty subset
of the set $~X~$.

$(ii)~$ Let $~Oxy~$ be a system of cartesian coordinates in a
plane. We consider the subsets $~Ox =~\{(x,0)\in {\bf
R}^{2}~|~(\forall)~x\in {\bf R}~\}~$ and $~Oy=~\{(0,y)\in {\bf
R}^{2}~|~(\forall)~y\in {\bf R}~\}~$ of $~X={\bf R}^{2}.~$ Let $~M
~=~\{ f_{1}=Id_{Ox}, f_{2}=Id_{Oy}, f_{3}=\sigma_{Ox} ,
f_{4}=\sigma_{Oy}~\}~\subset {\cal F}(S,{\bf R}^{2})~$ where the
functions $~f_{3} : Ox~\to ~Oy,~f_{4}~: Oy~\to~Ox~$ are defined by
$~f_{3}(x,0)=(0,x)~$ and $~f_{4}(0,y)=(y,0)~$ ($~\sigma_{Ox}~$
resp. $~\sigma_{Oy}~$ is called the {\it saltus function defined
on} $~x~$ {\it -axis~} resp. $~y~$-{\it axis} ).

The set of composable pairs of $~M~$ is\\
$~M_{(2)} = \{~(f_{1},f_{1}); (f_{1},f_{3});
(f_{2},f_{2});(f_{2},f_{4}); (
f_{3},f_{2});(f_{3},f_{4});(f_{4},f_{1});(f_{4},f_{3})~\}.~$

We have:\\
$\mu(f_{1},f_{1})= f_{1}\circ f_{1}=f_{1}\in M;~~~\mu(f_{1},f_{3})= f_{3}\circ f_{1}=f_{3}\in M;~~~\mu(f_{2},f_{2})= f_{2}\circ f_{2}=f_{2}\in M;~$\\
$\mu(f_{2},f_{4})= f_{4}\circ f_{2}=f_{4}\in M;~~~\mu(f_{3},f_{2})= f_{2}\circ f_{3}=f_{3}\in M;~~~\mu(f_{3},f_{4})= f_{4}\circ f_{3}=f_{1}\in M;~$\\
$\mu(f_{4},f_{1})= f_{1}\circ f_{4}=f_{4}\in
M;~~~\mu(f_{4},f_{3})= f_{3}\circ f_{4}=f_{2}\in M.~$

The unit set of $~M~$ is $~M_{0} =
\{~f_{1}=Id_{Ox},f_{2}=Id_{Oy}~\}.~$
The source and target map $~\alpha, \beta : M~\to~M_{0}~$ are given by\\
$\alpha(f_{1})=\beta(f_{1})=f_{1};~~\alpha(f_{2})=\beta(f_{2})=f_{2};~~\alpha(f_{3})=\beta(f_{4})=f_{1};~~\alpha(f_{4})=\beta(f_{3})=f_{2}.~$

We obtain that
$~(M=\{~f_{1},f_{2},f_{3},f_{4}~\},\alpha,\beta;M_{0}
=\{~f_{1},f_{2}~\}~)~$ is a monoidoid.\hfill$\b$

{\bf Definition 1.3.}~Let $~( G, \alpha, \beta, \mu;G_{0} )~$ be a
monoidoid endowed with an injective map
$~\iota:~G~\to~G,~x~\to~\iota(x),~$ called the {\it inversion map}
( we shall write $~x^{-1}~$ for $~\iota(x)~$).  We say that $~( G,
\alpha, \beta, \mu, \iota ; G_{0})~$ is a {\it groupoid}, if the
{\it inverses property} holds, i.e. for each $~x\in G~$ we have
$~(x^{-1},x),~(x,x^{-1})\in G_{(2)}~$ and $~x^{-1} x = \beta(x),~
x x^{-1} = \alpha(x).~$ \hfill$\b$

The subset $~G_{0} = \alpha(G) = \beta(G)~$ of $~G~$  is called
the {\it unit set} of $~G~$ and we say that $~G~$ is a $~G_{0}~$ -
{\it groupoid}. For all unit $~u\in G_{0}~$ we have $~\alpha(u) =
\beta(u) = \iota(u) = u.~$

A $~G_{0}~$-groupoid $~G~$ will be denoted by $~( G,\alpha,\beta;
G_{0})~$  or $~( G ; G_{0} ).~$ The maps $~\alpha, \beta, \mu ~$
and $~\iota~$ are called the {\it structure functions} of $~G.~$
The map
 $~(\alpha,\beta) :G~\to~G_{0}\times G_{0},~x~\to~(\alpha,\beta)(x) = ( \alpha(x),\beta(x))~$ is called the {\it anchor map} of the groupoid $~G.~$
A groupoid $~(G,\alpha,\beta;G_{0})~$ is called {\it transitive},
if the anchor map $~(\alpha,\beta) : G~\to~G_{0}\times G_{0}~$ is
surjective.

By {\it group bundle} we mean a $~G_{0}~$-groupoid $~G~$ such that
$~\alpha (x) = \beta (x)~$ for all $~x\in G~$. Moreover,a group
bundle is the union of its isotropy groups $~G(u) =
\alpha^{-1}(u), u \in G_{0}~$ (here, two elements may be composed
iff they lie in the same fiber $\alpha^{-1}(u)$~).

If $~( G,\alpha,\beta;G_{0})~$ is a groupoid then $~Is(G) = \{
x\in G ~|~ \alpha (x)=\beta (x) \}~$ is a group bundle, called the
{\it isotropy group bundle }  of  $~G.~$

{\bf Remark 1.1.} $~(i)~$ The definition of the Brandt groupoid is
essentially the same as the one given by A. Coste, P. Dazord and
A. Weinstein in [2].

$(ii)~$ A groupoid is a monoidoid in which every element is
invertible.\hfill$\b$

{\bf Example 1.3.}~$~(i)~$ We consider monoidoid $~M~= ~\{~e, 1,
0, -1~\}~$ over $~M_{0} = \{e,1~\},~$ see Example 1.1 (iv). This
monoidoid is not a groupoid, since the element $~0~$ is not
invertible.

$(ii)~$ Let $~M ~=~\{ f_{1}=Id_{Ox}, f_{2}=Id_{Oy}, f_{3} ,
f_{4}~\}~$ where $~Ox=~\{(x,0)\in {\bf R}^{2}~|~(\forall)~x\in
{\bf R}~\},~$ $~Oy=~\{(0,y)\in {\bf R}^{2}~|~(\forall)~y\in {\bf
R}~\},~$ $~f_{3} : Ox~\to
~Oy,~(x,0)~\linebreak\to~f_{3}(x,0)=(0,x)~$ and $~f_{4}~:
Oy~\to~Ox,~(0,y)~\to~f_{4}(0,y)=(y,0).~$ We have that
$~(M=\{~f_{1},f_{2},f_{3},f_{4}~\},\alpha,\beta;M_{0}
=\{~f_{1},f_{2}~\}~)~$ is a monoidoid, see Example 1.2 (ii).

We define the map $~\iota : M~\to~M~$ by taking
$~\iota(f_{1})=f_{1},~\iota(f_{2})=f_{2},~\iota(f_{3})=f_{4}~$ and
$~\iota(f_{4})=f_{3}.~$ It is easy to verify that $~(M,\alpha,
\beta, \mu, \iota; M_{0})~$ is a groupoid, called the {\it
groupoid of saltus functions defined on the axes of coordinates in
a plane}. We will denote this groupoid by $~{\cal F}_{(4;2)}({\bf
R}^{2}).~$

The groupoid $~(G={\cal F}_{(4;2)}({\bf R}^{2})=\{ f_{1}, f_{2},f_{3},f_{4}\},\alpha, \beta, \mu, \iota;G_{0}= \{ f_{1},f_{2}\})~$ is a transitive groupoid. Indeed,\\
for $~(f_{i},f_{i})\in G_{0}\times G_{0}~$ exists $~f_{i}\in G~$ such that $~\alpha(f_{i})=f_{i}~$ and $~\beta(f_{i})=f_{i}~$ for $~i=1,2;~$\\
for $~(f_{1},f_{2})\in G_{0}\times G_{0}~$ exists $~f_{3}\in G~$ such that $~\alpha(f_{3})=f_{1}~$ and $~\beta(f_{3})=f_{2},~~$ and\\
for $~(f_{2},f_{1})\in G_{0}\times G_{0}~$ exists $~f_{4}\in G~$
such that $~\alpha(f_{4})=f_{2}~$ and $~\beta(f_{4})=f_{1}.~$
Hence, the anchor map is surjective.\hfill$\b$

In the following proposition we summarize some properties of the
structure functions of a groupoid obtained directly from
definitions.

{\bf Proposition 1.2.}~{\it Let $~( G, \alpha, \beta, \mu, \iota ; G_{0})~$ be a groupoid and $~u\in G_{0}.~$ Then the following assertions hold:\\
$(1.5)~~~ \alpha (x^{-1}) = \beta (x), ~~ \beta ( x^{-1}) = \alpha (x)~$ and $~ (x^{-1})^{-1} = x,~$ for all $~x\in G;~$\\
$(1.6)~~~$({\it cancellation law} ) If $~x\cdot z_{1} = x\cdot z_{2}~$ ( resp., $~z_{1} \cdot x = z_{2}\cdot x~$ ), then $~z_{1} = z_{2};~$\\
$(1.7)~~~G(u)=\alpha^{-1}(u)\cap\beta^{-1}(u)=\{ x\in
G~|~\alpha(x)=\beta(x)=u~\}$ is a group under the restriction
of $\mu~$ to $~G(u),$ called the {\it isotropy group at} $u$ of $~G;$\\
$(1.8)~~$ If $~G~$ is transitive, then the isotropy groups
$~G(u),~u\in G_{0}~$ are isomorphes}.

{\bf Proof.}~Using the Proposition 1.1 ( each groupoid is a
monoidoid ) and the definitions, it is easy to prove that the
assertions $~(1.5)-(1.8)~$ are valid.\hfill$\b$

{\bf Definition 1.4.}~ Let $~( G, \alpha, \beta, \mu,\iota;
G_{0})~$ and $~( G^{\prime}, \alpha^{\prime}, \beta^{\prime},
\mu^{\prime}, \iota^{\prime}; G_{0}^{\prime} )~$ be two groupoids.
A {\it morphism of groupoids} or {\it groupoid morphism} from
$~G~$ into $~G^{\prime}~$ is a map $~f: G \longrightarrow
G^{\prime}~$  such that $~f(\mu(x,y)) = \mu^{\prime}(f(x),f(y))~$
for all $~(x,y)\in G_{(2)}.$ \hfill$\b$

A morphism of groupoids $~f :G~\to~G^{\prime}~$ such that the map
$~f~$ is bijective is called {\it isomorphism of groupoids} or
{\it groupoid isomorphism}.

The category $~{\cal BG}roid~$ of groupoids has as its objects all
groupoids $~(G,\alpha, \beta; G_{0})~$ and as morphisms from $~(G,
\alpha, \beta; G_{0})~$ to $~( G^{\prime}, \alpha^{\prime},
\beta^{\prime}; G_{0}^{\prime})~$ the set of morphisms of
groupoids.

{\bf Example 1.4.}~$~(i)~$ A group $~{\cal G}~$ having $~e~$ as
unit element is just a $~\{ e\}~$- groupoid and conversely, every
groupoid with one unit element is a group. It follows that the
category $~{\cal G}r~$ of groups is a subcategory of the category
$~{\cal BG}roid.~$

$~(ii)~$ Any nonempty set $~G_{0}~$ may be regarded as a nul
monoidoid on itself ( see, Example 1.1 (ii) ). If we consider the
map $~\iota: G_{0}~\to~G_{0}~$ as the identity map on $~G_{0},~$
we obtain that $~(G_{0},\alpha=\beta=Id_{G_{0}}, \mu,
\iota=Id_{G_{0}};G_{0})~$ is a groupoid, called the {\it nul
groupoid } associated to set $~G_{0}.~$

$(iii)~$ {\it Disjoint union of two groupoids}.~ Let
$~(G_{i},\alpha_{i}, \beta_{i}, \mu_{i}, \iota_{i}; G_{0,i})~, i =
1,2~$ be two groupoids such that $~G_{1}\cap G_{2} = \emptyset.~$

We consider the disjoint union $ ~( G = G_{1}\coprod G_{2},
\alpha, \beta,\mu ;G_{0}=G_{1,0}\coprod G_{2,0})~$ of the
monoidoids $~G_{1}~$ and $~G_{2}~$ ( see, Example 1.1 (iii) ).

We define the map $~\iota : G~\to~G~$ by taking $~\iota(x) =
\iota_{1}(x),~$ if $~x\in G_{1}~$ and $~\iota(x) = \iota_{2}(x)~$
if $~x\in G_{2}.~$ We have that $ ~(G = G_{1}\coprod G_{2},
\alpha, \beta,\mu ,\iota ;G_{0}=G_{1,0}\coprod G_{2,0})~$ is a
groupoid, called the {\it disjoint union } of the groupoids
$~G_{1}~$ and $~G_{2}.~$

In particular, the disjoint union of groups $~G_{i},~ i=1,2~$ is a
groupoid, called the {\it groupoid associated to groups} $~G_{i},~
i=1,2~$ ( for this groupoid, the unit set
 is $~G_{0}= \{ e_{1}, e_{2}\}~$ where $~e_{i}~$ is the unity of $~G_{i}, i = 1,2~$ ).\hfill$\b$

A finite groupoid $~(G;G_{0})~$ such that $~|~G~| = n~$ and
$~|~G_{0}~| =m~$ is called {\it $~( n ; m )-~$ groupoid} or {\it
finite groupoid of type $~(n;m)~$}.

{\bf Example 1.5.}~$~(i)~$ Each finite groupoid of type $~(n;1)~$
is a group.

$(ii)~$ Each finite groupoid of type $~(n;n)~$ is a nul groupoid.

$(iii)~$ The groupoid $~{\cal F}_{(4;2)}({\bf R}^{2})~$ ( see,
Example 1.3 (ii) ) is a finite transitive groupoid of type
$~(4;2).~$

$(iv)~$ Let be the Klein $~4~$-group $~K_{4} = \{~(1),~\sigma =
(12)(34),~\tau=(13)(24),~\sigma\circ \tau = (14)(23)~\}\subset
S_{4}~$  ( it is a subgroup of the symmetric group $~S_{4}~$ of
degree $~4~$). We have $~\sigma^{2} = \tau^{2} = (1)~$ and
$~\tau\circ \sigma =\sigma\circ \tau.$

We consider the disjoint union $~G = K_{4}\coprod {\bf Z}_{4}~$ of
the groups $~K_{4}~$ and $~{\bf Z}_{4},~$ where $~{\bf Z}_{4}=
\{~\widehat{0}, \widehat{1},\widehat{2}, \widehat(3)~\}~$ is the
group of congruences classes of integers modulo $~4.~$ We  obtain
a groupoid $~G~$ of order $~8~$ with unit set $~G_{0} = \{~(1),
\widehat{0}~\}~$ of type $~(8;2).~$

The structure functions $~\alpha,~\beta~, \iota~$ and $~\mu~$  of $~G = K_{4}\coprod {\bf Z}_{4}~$ are given in the following tables:\\
$$\begin{array}{|c|c|c|c|c|c|c|c|c|c|} \hline
x         & (1) & \sigma & \tau  & \sigma\circ \tau & \widehat{0}
& \widehat{1} & \widehat{2} & \widehat{3} \\ \hline \alpha(x) &
(1) & (1)    &  (1)  & (1)              & \widehat{0} &
\widehat{0} & \widehat{0}  & \widehat{0} \cr \hline \beta(x)  &
(1) & (1)    &  (1)  & (1)              & \widehat{0} &
\widehat{0} & \widehat{0}  & \widehat{0} \cr \hline \iota(x)  &
(1) & \sigma & \tau  & \sigma\circ\tau  & \widehat{0} &
\widehat{3} & \widehat{2} &  \widehat{1} \cr \hline
\end{array}$$ $$\begin{array}{|c|c|c|c|c|c|c|c|c|c|} \hline
\mu               & (1)              & \sigma          & \tau
& \sigma\circ \tau & \widehat{0} & \widehat{1} & \widehat{2} &
\widehat{3} \\ \hline (1)               & (1)              &
\sigma          & \tau              & \sigma\circ \tau &
&             &             &             \cr \hline \sigma
& \sigma           & (1)             & \tau\circ\sigma   & \tau
&             &             &             &             \cr \hline
\tau              & \tau             & \sigma\circ\tau &  (1)
&  \sigma          &             &             &             &
\cr \hline \sigma\circ\tau   & \sigma\circ\tau  & \tau
& \sigma            &  (1)             &             &
&             &             \cr \hline \widehat{0}       &
&                 &                   &                  &
\widehat{0} & \widehat{1} & \widehat{2} & \widehat{3} \cr \hline
\widehat{1}       &                  &                 &
&                  & \widehat{1} & \widehat{2} & \widehat{3} &
\widehat{0} \cr \hline \widehat{2}       &                  &
&                   &                  & \widehat{2} & \widehat{3}
& \widehat{0} & \widehat{1} \cr \hline \widehat{3}       &
&                 &                   &                  &
\widehat{3} & \widehat{0} & \widehat{1} & \widehat{2} \cr \hline
\end{array}$$\hfill$\b$\\[0.1cm]

{\bf 2. The program $~BGroidAP1~$ for to test if an universal algebra is a groupoid}\\[0.1cm]

We consider a given finite universal algebra $~(G, \alpha, \beta,
\mu , \iota;G_{0})~$ such that $~|~G~|~ = n~$ and $~|~G_{0}~|~ =
m~$ with $~1\leq m\leq n.~$ We denote the elements of $~G~$ by
$~a_{1}, a_{2}, \cdots, a_{m}, a_{m+1}, \cdots, a_{n}~$ such that
$~G_{0} =\{~a_{1}, a_{2}, \cdots, a_{m}~\}.~$ Hence, the elements
$~a_{k},~ k = \overline{1,m}~$ are the units of $~G.~$

We give an algorithm for decide if the universal algebra $~(
G,\alpha, \beta, \mu, \iota; G_{0})~$ is a $~G_{0}~$- groupoid.
This algorithm is constituted by the following stages.

{\bf Stage I.}~{\it We introduce the initial data}: $~~~n~$- the
number of elements of $G;~~~m~$- the number of elements of
$G_{0};~$ the functions $~\alpha, \beta, \iota~$ and $~\mu~$ given
by its tables of structure.

{\bf Stage II.}~{\it Test if the universal algebra $~(G, \alpha,
\beta, \mu, \iota; G_{0})~$ considered
in the first stage is a groupoid}. For this, the following steps are executed:\\
{\bf step 1.} $ (G,\alpha, \beta, \mu,\iota )~$ is a structure
well-defined, i.e. the functions $~\alpha,~\beta~$ are
surjections, $~\iota~$ is injective and
$~\mu~$ is defined on the composable pairs $~G_{(2)}~$ with values in $~G;~$\\
{\bf step 2.} $(G,\alpha, \beta, \mu )~$ is a semigroupoid, i.e. the multiplication law $~\mu~$ is associative;\\
{\bf step 3.}  the semigroupoid $~(G; G_{0})~$ is a monoidoid,
i.e. the identities property
 is verified;\\
{\bf step 4.} the monoidoid $~(G; G_{0} )~$ is a groupoid, i.e. each element of $~G~$ is invertible.\\
{\bf step 5.} If the above steps are satisfied, make the tables of
the structure functions $~\alpha, \beta, \iota~$ and $~\mu~$
and write the mesage "{\it $G$ is a groupoid}$\,$".\\

Let us we present the correspondence between the initial data and
input data:
\begin{align*}
G=\{~a_{1},a_{2},\ldots,\ldots,a_{m},a_{m+1},\ldots,
a_{n}~\}&\longleftrightarrow\{~1,2,\ldots,m,m+1,\ldots,n~\}\\
Initial\ data &\longleftrightarrow Input\ data\\
|~G~|=n&\longleftrightarrow  n\\
|~G_{0}~|=m&\longleftrightarrow m\\
\begin{array}{|c|c|c|c|c|c|c|}\hline
 a_{k}    & \!a_{1}\!\!\! & \!\!\cdots\!\! & \!a_{m}\!\! & a_{m+1} & \!\!\cdots\!\! & a_{n} \\ \hline
  \!\alpha(a_{k})\!\!& \!a_{1}\!\!\! & \!\!\cdots\!\! & \!a_{m}\!\! & \alpha(a_{m+1}) & \!\!\cdots\!\! &
\!\alpha(a_{n})\!\! \cr \hline
 \! \beta(a_{k})\!\!&\! a_{1}\!\!\! &\!\! \cdots\!\! & \!a_{m}\!\! & \beta(a_{m+1}) & \!\!\cdots\!\! & \beta(a_{n}) \cr \hline
\!\iota(a_{k})\!\! & \!a_{1}\!\!\! & \!\!\cdots\!\! & \!a_{m}\!\!
& \iota(a_{m+1}) & \!\!\cdots \!\!& \iota(a_{n}) \cr \hline
\end{array}&\longleftrightarrow \begin{array}{cccccc} \\
       1 \!\!&\!\! \cdots\!\! & m & \!\!u_{l}(m+1)\!\! & \!\!\cdots\!\! & \!\!u_{l}(n)\cr
       1 \!\!&\!\! \cdots\!\! & m & \!\!u_{r}(m+1)\!\! & \!\!\cdots\!\! & \!\!u_{r}(n) \cr
       1 \!\!&\!\! \cdots\!\! & m & \!\!inv(m+1) \!\!& \!\!\cdots\!\! & \!\!inv(n)\cr
\end{array}\\
\begin{array}{|r|c|c|c|c|c|}\hline
\mu  & a_{1} & \cdots & a_{k}            & \cdots & a_{n} \\
\hline a_{1}  &       &        &                  &        &
\cr \hline \cdots &       &        &                  &        &
\cr \hline a_{j}  &       &        & a_{j}\cdot a_{k} &        &
\cr \hline \cdots &       &        &                  &        &
\cr \hline a_{n}  &       &        &                  &        &
\cr \hline
\end{array}&\longleftrightarrow \begin{array}{cccccc} \\
a_{11}  & \cdots    &  a_{1k}   & \cdots  & a_{1n} \cr a_{21}  &
\cdots    & a_{2k}    & \cdots  &  a_{2n}\cr \cdots  & \cdots    &
\cdots    & \cdots  & \cdots  \cr a_{j1}  &  \cdots   & a_{jk} &
\cdots  & a_{jn} \cr \cdots  & \cdots    & \cdots    & \cdots  &
\cdots  \cr a_{n1}  & \cdots    & a_{nk}    & \cdots  & a_{nn}.
\cr
\end{array}
\end{align*}

In the table of $~\mu~$ the element$~\mu(a_{j},
a_{k})~=~a_{j}\cdot a_{k}~$ is defined iff $~\beta(a_{j}) =
\alpha(a_{k}).~$ The absence of an element from the arrow $~"j"~$
and the column $~"k"~$ indicates the fact that the pair $~(a_{j},
a_{k})\in G\times G~$ is not composable. The element $~a_{jk}~$ is
represented by $~0~$ if the product $~a_{j}\cdot a_{k}~$ is not
defined.

{\bf Example 2.1.}~ Consider the groupoid $~G = {\cal
F}_{(4;2)}{\bf R}^{2}~$ given in Example 1.3 (ii). We have
\begin{align*}
G=\{~f_{1}=Id_{Ox},f_{2}=Id_{Oy},f_{3}=\sigma_{Ox},f_{4}=\sigma_{Oy}~\}&\longleftrightarrow\{~1,2,3,4~\}\\
Initial\ data &\longleftrightarrow Input\ data\\
|~G~|=4&\longleftrightarrow 4\\
|~G_{0}~|=2&\longleftrightarrow 2\end{align*}
\begin{align*}
\begin{array}{|c|c|c|c|c|}\hline
 f_{k}    & f_{1} & f_{2} & f_{3} & f_{4} \\ \hline
  \alpha(f_{k})& f_{1} & f_{2} & f_{1} & f_{2}  \cr \hline
  \beta(f_{k})& f_{1} & f_{2} & f_{2} & f_{1} \cr \hline
\iota(f_{k}) & f_{1} & f_{2} & f_{4} & f_{3} \cr \hline
\end{array}&\longleftrightarrow\begin{array}{cccccc} \\
       1 & 2 & 1 & 2\cr
       1 & 2 & 2 & 1 \cr
       1 & 2 & 4 & 3\cr
\end{array}\\
\begin{array}{|r|c|c|c|c|}\hline
\mu    & f_{1}  &  f_{2} & f_{3} & f_{4} \\ \hline f_{1}  & f_{1}
&        & f_{3} &         \cr \hline f_{2}  &        &  f_{2} &
&  f_{4}  \cr \hline f_{3}  &        &  f_{3} &       &  f_{1}
\cr \hline f_{4}  & f_{4}  &        & f_{2} &         \cr \hline
\end{array}&\longleftrightarrow\begin{array}{cccc} \\
1  &  0   &  3   &  0 \cr 0  &  2   &  0   &  4 \cr 0  &  3   &  0
&  1 \cr 4  &  0   &  2   &  0  \cr
\end{array}\end{align*} \hfill$\b$\\[0.2cm]

The implementation of the above algorithm on computer is realized
in the program $~BGroidAP1.~$ This program is composed from two
modules denoted by $~unit11.dfm~$ and $~unit11.pas~$. The module
$~unit11.pas~$ is composed from the principal program
followed of procedures and functions. This is presented as follows.\\

unit Unit1;\\[-0.2cm]

interface\\[-0.2cm]

uses

\hspace*{0.5cm}Windows, Messages, SysUtils, Classes, Graphics,
Controls, Forms, Dialogs,

\hspace*{0.5cm}Grids, DBGrids, ShellAPI, Db, DBTables, StdCtrls,
Menus, ExtCtrls,

\hspace*{0.5cm}ComCtrls, ToolWin, Spin;\\[-0.2cm]

const

\hspace*{1.3cm}nmax = 200;\\[-0.2cm]

type

\hspace*{0.5cm}TForm1 = class(TForm)

\hspace*{0.8cm}MainMenu1: TMainMenu;

\hspace*{0.8cm}File1: TMenuItem;

\hspace*{0.8cm}OpenFile1: TMenuItem;

\hspace*{0.8cm}SaveFile1: TMenuItem;

\hspace*{0.8cm}GroupBox1: TGroupBox;

\hspace*{0.8cm}StringGrid1: TStringGrid;

\hspace*{0.8cm}StringGrid2: TStringGrid;

\hspace*{0.8cm}OpenDialog1: TOpenDialog;

\hspace*{0.8cm}SaveDialog1: TSaveDialog;

\hspace*{0.8cm}ToolBar1: TToolBar;

\hspace*{0.8cm}Splitter4: TSplitter;

\hspace*{0.8cm}ToolButton1: TToolButton;

\hspace*{0.8cm}New1: TMenuItem;

\hspace*{0.8cm}ToolBar3: TToolBar;

\hspace*{0.8cm}ToolButton4: TToolButton;

\hspace*{0.8cm}ToolButton5: TToolButton;

\hspace*{0.8cm}Label2: TLabel;

\hspace*{0.8cm}SpinEdit1: TSpinEdit;

\hspace*{0.8cm}ToolButton6: TToolButton;

\hspace*{0.8cm}Label3: TLabel;

\hspace*{0.8cm}SpinEdit2: TSpinEdit;

\hspace*{0.8cm}StatusBar1: TStatusBar;

\hspace*{0.8cm}ToolButton3: TToolButton;

\hspace*{0.8cm}ToolButton11: TToolButton;

\hspace*{0.8cm}procedure FormShow(Sender: TObject);

\hspace*{0.8cm}procedure Button2Click(Sender: TObject);

\hspace*{0.8cm}procedure StringGrid1SetEditText(Sender: TObject;
ACol, ARow: Integer;

\hspace*{1.1cm}const Value: String);

\hspace*{0.8cm}procedure StringGrid2SetEditText(Sender: TObject;
ACol, ARow: Integer;

\hspace*{1.1cm}const Value: String);

\hspace*{0.8cm}procedure OpenFile1Click(Sender: TObject);

\hspace*{0.8cm}procedure SaveFile1Click(Sender: TObject);

\hspace*{0.8cm}procedure New1Click(Sender: TObject);

\hspace*{0.8cm}procedure ToolButton4Click(Sender: TObject);

\hspace*{0.8cm}procedure ToolButton3Click(Sender: TObject);

\hspace*{0.5cm}private

\hspace*{0.8cm}err\_message : String;

\hspace*{1.4cm}m, n : Integer;

\hspace*{0.8cm}h : array[0..nmax, 0..nmax] of Byte;

\hspace*{0.8cm}u\_left, u\_right, inv : array[0..nmax] of Integer;\\[-0.2cm]

\hspace*{0.8cm}procedure WMDropFiles(var Msg: TWMDropFiles);
message WM\_DROPFILES;

\hspace*{0.8cm}procedure PerformFileOpen(const FileName1 :
string);

\hspace*{0.8cm}procedure PerformFileSave(const FileName1 :
string);

\hspace*{0.8cm}procedure MakeUnitsTable;

\hspace*{1.4cm}\qquad procedure MakeGroupoidTable;

\hspace*{0.8cm}function ToStr(x : Integer) : String;

\hspace*{0.8cm}function IsStructure : Boolean;

\hspace*{0.8cm}function IsSemigroupoid : Boolean;

\hspace*{0.8cm}function IsMonoidoid : Boolean;

\hspace*{0.8cm}function IsGroupoid : Boolean;

\hspace*{0.5cm}public\\[-0.2cm]

\hspace*{0.5cm}end;\\[-0.2cm]

var

\hspace*{0.5cm}Form1: TForm1;\\[-0.2cm]

implementation\\[-0.2cm]

\{\$R *.DFM\}\\[-0.2cm]

procedure TForm1.FormShow(Sender: TObject);

var

\hspace*{1.4cm}i, j : Byte;

begin

\hspace*{0.5cm}DragAcceptFiles(Handle, True);

\hspace*{0.5cm}StringGrid1.EditorMode := True;

\hspace*{0.5cm}n := 0;

\hspace*{0.5cm}m := 0;

\hspace*{0.5cm}for i := 0 to nmax do

\hspace*{1.4cm}for j := 0 to nmax do

\hspace*{1.4cm}h[i, j] := 0;

\hspace*{0.5cm}for i := 0 to nmax do begin

\hspace*{0.8cm}u\_left[i] := 0;

\hspace*{0.8cm}u\_right[i] := 0;

\hspace*{0.8cm}inv[i] := 0;

\hspace*{0.5cm}end;

end;\\[-0.2cm]

function TForm1.ToStr;

var

\hspace*{1.4cm}ss : String;

begin

\hspace*{0.5cm}str(x, ss);

\hspace*{0.5cm}ToStr := ss;

end;\\[-0.2cm]

procedure TForm1.MakeUnitsTable;

var

\hspace*{1.4cm}i : Integer;

begin

\hspace*{1.4cm}StringGrid2.RowCount := 4;

\hspace*{1.4cm}StringGrid2.ColCount := n + 1;

\hspace*{1.4cm}StringGrid2.Cells[0, 1] := 'u\_l';

\hspace*{1.4cm}StringGrid2.Cells[0, 2] := 'u\_r';

\hspace*{1.4cm}StringGrid2.Cells[0, 3] := 'inv';

\hspace*{0.5cm}for i := 1 to n do begin

\hspace*{1.4cm}StringGrid2.Cells[i, 0] := tostr(i);

\hspace*{1.4cm}StringGrid2.Cells[i, 1] := tostr(u\_left[i]);

\hspace*{1.4cm}StringGrid2.Cells[i, 2] := tostr(u\_right[i]);

\hspace*{1.4cm}StringGrid2.Cells[i, 3] := tostr(inv[i])

\hspace*{0.5cm}end

end;\\[-0.2cm]

procedure TForm1.MakeGroupoidTable;

var

\hspace*{0.5cm}i, j : Integer;

begin

\hspace*{0.5cm}StatusBar1.SimpleText := 'G has not been tested';

\hspace*{0.5cm}if n $ > $ 0 then begin

\hspace*{0.8cm}StringGrid1.RowCount := n + 1;

\hspace*{0.8cm}StringGrid1.ColCount := n + 1;

\hspace*{0.8cm}for i := 1 to n do begin

\hspace*{1cm}StringGrid1.Cells[0, i] := tostr(i);

\hspace*{1cm}StringGrid1.Cells[i, 0] := tostr(i);

\hspace*{0.8cm}end;

\hspace*{0.8cm}for i := 1 to n do

\hspace*{1cm}for j := 1 to n do

\hspace*{1.4cm}if h[i, j] $ < > $  0 then

\hspace*{1.7cm}StringGrid1.Cells[j, i] := tostr(h[i, j])

\hspace*{1.4cm}else

\hspace*{1.7cm}StringGrid1.Cells[j, i] := ''

\hspace*{0.5cm}end else begin

\hspace*{0.8cm}StringGrid1.RowCount := 2;

\hspace*{0.8cm}StringGrid1.ColCount := 2

\hspace*{0.5cm}end;

\hspace*{0.5cm}MakeUnitsTable

end;\\[-0.2cm]

function TForm1.IsStructure;

var

\hspace*{1.4cm}i, j : Byte;

begin

\hspace*{1.4cm}IsStructure := true;

\hspace*{0.5cm}for i := 1 to n do begin

\hspace*{1.4cm}if (u\_left[i] = 0) or (u\_right[i] = 0) or (inv[i]
= 0) then begin

\hspace*{3.1cm}IsStructure := false;

\hspace*{1cm}err\_message := 'Structure incomplete';

\hspace*{1cm}exit

\hspace*{0.8cm}end

\hspace*{0.5cm}end;

\hspace*{0.5cm}for i := 1 to n do for j := 1 to n do

\hspace*{1.4cm}if (u\_right[i] = u\_left[j]) and (h[i, j] = 0)
then begin

\hspace*{3.1cm}IsStructure := false;

\hspace*{1cm}err\_message := 'Structure incomplete';

\hspace*{1cm}exit

\hspace*{0.8cm}end

end;\\[-0.2cm]

function TForm1.IsSemigroupoid;

var

\hspace*{1.4cm}i, j, k : Byte;

begin

\hspace*{1.4cm}if IsStructure then begin

\hspace*{0.8cm}IsSemigroupoid := true;

\hspace*{0.8cm}for i := 1 to n do for j := 1 to n do if
u\_right[i] = u\_left[j] then

\hspace*{1cm}for k := 1 to n do if u\_right[j] = u\_left[k] then

\hspace*{1.4cm}if h[h[i, j], k] $ < > $ h[i, h[j, k]] then begin

\hspace*{1.7cm}IsSemigroupoid := false;

\hspace*{1.7cm}err\_message := tostr(i) + ', ' + tostr(j) + ', ' +
tostr(k) + ' - not asociative';

\hspace*{1.7cm}exit

\hspace*{1.4cm}end

\hspace*{0.5cm}end else IsSemigroupoid := false

end;\\[-0.2cm]

function TForm1.IsMonoidoid;

var

\hspace*{1.4cm}i : Byte;

begin

\hspace*{1.4cm}if IsSemigroupoid then begin

\hspace*{0.8cm}IsMonoidoid := true;

\hspace*{0.8cm}for i := 1 to n do

\hspace*{1cm}if (h[u\_left[i], i] $ < > $ i) or (h[i, u\_right[i]]
$ < > $ i) then begin

\hspace*{1.4cm}IsMonoidoid := false;

\hspace*{1.4cm}err\_message := tostr(i) + ' has no unit';

\hspace*{1.4cm}exit

\hspace*{1cm}end

\hspace*{0.5cm}end else IsMonoidoid := false

end;\\[-0.2cm]

function TForm1.IsGroupoid;

var

\hspace*{1.4cm}i : Byte;

begin

\hspace*{1.4cm}if IsMonoidoid then begin

\hspace*{0.8cm}IsGroupoid := true;

\hspace*{0.8cm}for i := 1 to n do begin

\hspace*{1cm}if u\_right[i] = u\_left[inv[i]] then

\hspace*{1.4cm}if h[i, inv[i]] $ < > $ u\_left[i] then begin

\hspace*{1.7cm}IsGroupoid := false;

\hspace*{1.7cm}err\_message := tostr(i) + ' has no inverse';

\hspace*{1.7cm}exit

\hspace*{1.4cm}end;

\hspace*{1cm}if u\_right[inv[i]] = u\_left[i] then

\hspace*{1.4cm}if h[inv[i], i] $ < > $ u\_right[i] then begin

\hspace*{1.7cm}IsGroupoid := false;

\hspace*{1.7cm}err\_message := tostr(i) + ' - has no inverse';

\hspace*{1.7cm}exit

\hspace*{1.4cm}end

\hspace*{0.8cm}end

\hspace*{0.5cm}end else IsGroupoid := false

end;\\[-0.2cm]

procedure TForm1.WMDropFiles(var Msg: TWMDropFiles);

var

\hspace*{0.5cm}CFileName: array[0..MAX\_PATH] of Char;

begin

\hspace*{0.5cm}try

\hspace*{0.8cm}if DragQueryFile(Msg.Drop, 0, CFileName, MAX\_PATH)
$ > $ 0 then

\hspace*{0.8cm}begin

\hspace*{1cm}\{CheckFileSave;\}

\hspace*{1cm}PerformFileOpen(CFileName);

\hspace*{1cm}Msg.Result := 0;

\hspace*{0.8cm}end;

\hspace*{0.5cm}finally

\hspace*{0.8cm}DragFinish(Msg.Drop);

\hspace*{0.5cm}end;

end;\\[-0.2cm]

procedure TForm1.PerformFileOpen(const FileName1 : string);

var

\hspace*{1.4cm}fin : TextFile;

\hspace*{0.5cm}i, j : Integer;

begin

\hspace*{1.4cm}AssignFile(fin, FileName1);

\hspace*{0.5cm}reset(fin);

\hspace*{0.5cm}readln(fin, n);

\hspace*{0.5cm}readln(fin, m);

\hspace*{0.5cm}readln(fin);

\hspace*{0.5cm}for i := 1 to n do

\hspace*{1.4cm}read(fin, u\_left[i]);

\hspace*{0.5cm}for i := 1 to n do

\hspace*{1.4cm}read(fin, u\_right[i]);

\hspace*{0.5cm}for i := 1 to n do

\hspace*{1.4cm}read(fin, inv[i]);

\hspace*{0.5cm}for i := 1 to n do

\hspace*{1.4cm}for j := 1 to n do

\hspace*{1.4cm}read(fin, h[i, j]);

\hspace*{0.5cm}CloseFile(fin);

\hspace*{0.5cm}MakeGroupoidTable;

end;\\[-0.2cm]

procedure TForm1.PerformFileSave(const FileName1 : string);

var

\hspace*{1.4cm}f : TextFile;

\hspace*{0.5cm}i, j : Integer;

begin

\hspace*{1.4cm}AssignFile(f, FileName1);

\hspace*{0.5cm}rewrite(f);

\hspace*{0.5cm}writeln(f, n);

\hspace*{0.5cm}writeln(f, m);

\hspace*{0.5cm}writeln(f);

\hspace*{0.5cm}for i := 1 to n do

\hspace*{1.4cm}write(f, u\_left[i], ' ');

\hspace*{0.5cm}writeln(f);

\hspace*{0.5cm}for i := 1 to n do

\hspace*{1.4cm}write(f, u\_right[i], ' ');

\hspace*{0.5cm}writeln(f);

\hspace*{0.5cm}for i := 1 to n do

\hspace*{1.4cm}write(f, inv[i], ' ');

\hspace*{0.5cm}writeln(f);

\hspace*{0.5cm}writeln(f);

\hspace*{0.5cm}for i := 1 to n do begin

\hspace*{1.4cm}for j := 1 to n do

\hspace*{1.4cm}write(f, h[i, j], ' ');

\hspace*{0.8cm}writeln(f);

\hspace*{0.5cm}end;

\hspace*{0.5cm}CloseFile(f)

end;\\[-0.2cm]

procedure TForm1.Button2Click(Sender: TObject);

begin

\hspace*{1.4cm}MakeGroupoidTable;

\hspace*{1.4cm}err\_message := '';

\hspace*{1.4cm}if IsGroupoid then

\hspace*{1.4cm}StatusBar1.SimpleText := 'G is a groupoid'

\hspace*{0.5cm}else

\hspace*{1.4cm}StatusBar1.SimpleText := err\_message;

end;\\[-0.2cm]

procedure TForm1.StringGrid1SetEditText(Sender: TObject; ACol,

\hspace*{0.5cm}ARow: Integer; const Value: String);

var

\hspace*{1.4cm}nr : Byte;

\hspace*{0.5cm}cod : Integer;

begin

\hspace*{1.4cm}Val(Value, nr, cod);

\hspace*{0.5cm}if cod $ < > $ 0 then

\hspace*{1.4cm}h[ARow, ACol] := 0

\hspace*{0.5cm}else

\hspace*{3.1cm}h[ARow, ACol] := nr

end;\\[-0.2cm]

procedure TForm1.StringGrid2SetEditText(Sender: TObject; ACol,

\hspace*{0.5cm}ARow: Integer; const Value: String);

var

\hspace*{1.4cm}nr : Byte;

\hspace*{0.5cm}cod : Integer;

begin

\hspace*{0.5cm}Val(Value, nr, cod);

\hspace*{0.5cm}if cod $ < > $ 0 then nr := 0;

\hspace*{1.4cm}case ARow of

\hspace*{1.4cm}1 : u\_left[ACol] := nr;

\hspace*{1.4cm}2 : u\_right[ACol] := nr;

\hspace*{1.4cm}3 : inv[ACol] := nr

\hspace*{0.5cm}end

end;\\[-0.2cm]

procedure TForm1.OpenFile1Click(Sender: TObject);

begin

\hspace*{1.4cm}if OpenDialog1.Execute then

\hspace*{3.1cm}PerformFileOpen(OpenDialog1.FileName)

end;\\[-0.2cm]

procedure TForm1.SaveFile1Click(Sender: TObject);

begin

\hspace*{1.4cm}if SaveDialog1.Execute then

\hspace*{1.4cm}PerformFileSave(SaveDialog1.FileName)

end;\\[-0.2cm]

procedure TForm1.New1Click(Sender: TObject);

begin

\hspace*{0.5cm}SpinEdit1.Value := n;

\hspace*{0.5cm}SpinEdit2.Value := m;

\hspace*{1.4cm}ToolBar3.Show;

\hspace*{0.5cm}ToolBar1.Hide

end;\\[-0.2cm]

procedure TForm1.ToolButton3Click(Sender: TObject);

begin

\hspace*{1.4cm}ToolBar1.Show;

\hspace*{1.4cm}ToolBar3.Hide

end;\\[-0.2cm]

procedure TForm1.ToolButton4Click(Sender: TObject);

var

\hspace*{1.4cm}i, j : Byte;

begin

\hspace*{1.4cm}n := SpinEdit1.Value;

\hspace*{1.4cm}m := SpinEdit2.Value;

\hspace*{0.5cm}for i := 1 to n do

\hspace*{1.4cm}for j := 1 to n do

\hspace*{1.4cm}h[i, j] := 0;

\hspace*{0.5cm}ToolBar1.Show;

\hspace*{0.5cm}ToolBar3.Hide;

\hspace*{0.5cm}MakeGroupoidTable

end;

end.\\[0.1cm]

We illustrate the utilisation of the program $~BGgroidAP1~$ in the following examples.\\

{\bf Example 2.2.}~ Consider the universal algebra $~( G, \alpha,
\beta, \mu, \iota; G_{0}),~$
where \\
$G =\{~x_{1}, x_{2}, x_{3}, x_{4}, x_{5} , x_{6} , x_{7},x_{8},
x_{9}~\},~$ $~G_{0} = \{~ x_{1}, x_{2}, x_{3}~\}~$ and the inputs
data are the following:
$$\begin{array}{ccccccccc} \\
 9 & & & & & & & &\\
 3 & & & & & & & &\\[0.1cm]

1 & 2 & 3 & 1 & 1 & 2 & 2 & 3 & 3 \\
1 & 2 & 3 & 2 & 3 & 1 & 3 & 1 & 2\\
1 & 2 & 3 & 6 & 8 & 4 & 9 & 5 & 7\\
\end{array}~~~~~~~~~~\begin{array}{ccccccccc} \\
 1 & 0 & 0 & 4 & 5 & 0 & 0 & 0 & 0 \\
 0 & 2 & 0 & 0 & 0 & 6 & 7 & 0 & 0 \\
 0 & 0 & 3 & 0 & 0 & 0 & 0 & 8 & 9 \\
 0 & 4 & 0 & 0 & 0 & 1 & 5 & 0 & 0 \\
 0 & 0 & 5 & 0 & 0 & 0 & 0 & 1 & 4 \\
 6 & 0 & 0 & 2 & 7 & 0 & 0 & 0 & 0 \\
 0 & 0 & 7 & 0 & 0 & 0 & 0 & 6 & 2 \\
 8 & 0 & 0 & 9 & 3 & 0 & 0 & 0 & 0 \\
 0 & 9 & 0 & 0 & 0 & 8 & 3 & 0 & 0 \\
\end{array}$$

Execute the program $~BGroidAP1~$ for
$~G=\{~x_{j}~|~j=\overline{1,9}~\}~$ and the window program of
obtained results is presented in the Figure \ref{fig1}.

\begin{figure}[!htb]
  \centering
  \includegraphics[scale=1]{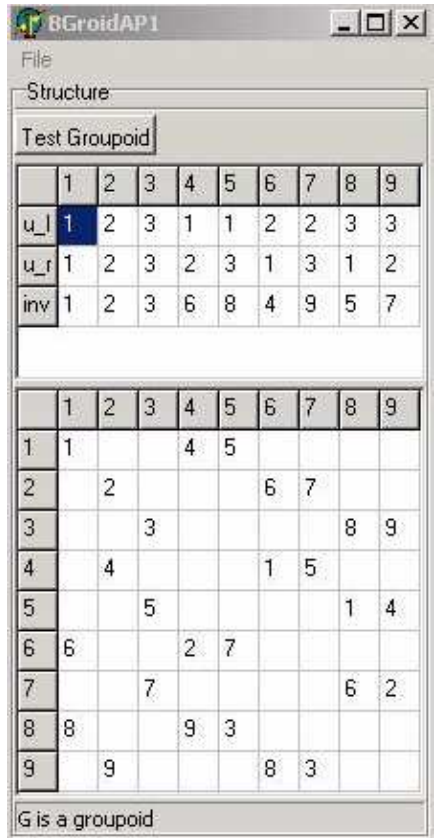}
  \caption{A $(9;3)$ - groupoid} \label{fig1}
\end{figure}

Therefore, $~G~$ is a groupoid of type $~(9;3).~$\hfill$\b$

{\bf Example 2.3.}~ Consider the universal algebra $~( G, \alpha,
\beta, \mu, \iota; G_{0}),~$ where $~G =\{~a_{1}, a_{2}, a_{3},
a_{4}, a_{5} , a_{6} ~\},~$ $~G_{0} = \{~ a_{1}~\}~$ and the
inputs data are the following:
$$\begin{array}{cccccc} \\
 6 & & & & & \\
 1 & & & & & \\[0.1cm]

1 & 1 & 1 & 1 & 1 & 1 \\
1 & 1 & 1 & 1 & 1 & 1 \\
1 & 3 & 2 & 4 & 6 & 5 \\
\end{array}~~~~~~~~~~~~~~~~~\begin{array}{cccccc} \\
 1 & 2 & 3 & 4 & 5 & 6 \\
 2 & 3 & 1 & 6 & 4 & 5 \\
 3 & 1 & 2 & 5 & 6 & 4 \\
 4 & 5 & 6 & 1 & 2 & 3 \\
 5 & 6 & 4 & 3 & 1 & 2 \\
 6 & 4 & 5 & 2 & 3 & 1 \\
\end{array}$$ \hfill$\b$

Execute the program $~BGroidAP1~$ for
$~G=\{~a_{k}~|~k=\overline{1,6}~\}~$ and the window program of
obtained results is presented in the Figure \ref{fig2}.

\begin{figure}[!htb]
  \centering
  \includegraphics[scale=1]{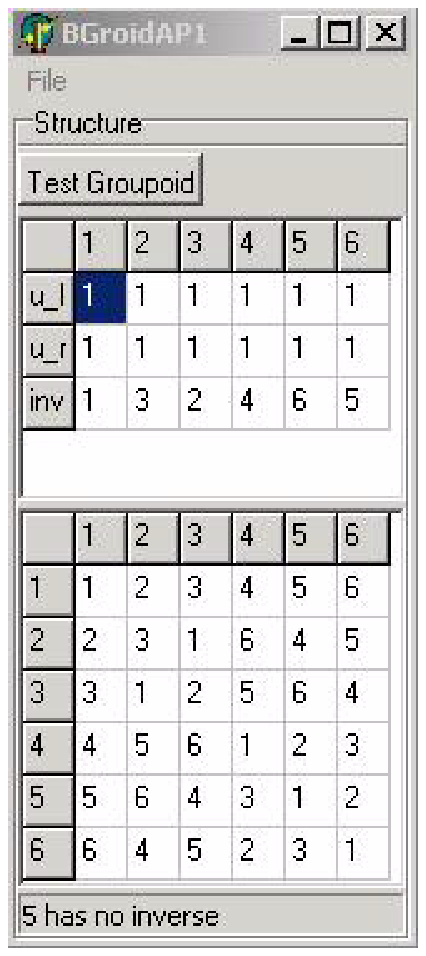}
  \caption{A monoid which is not a group} \label{fig2}
\end{figure}

Therefore, $~G~$ is not a groupoid of type $~(6;1)~$ ( hence ,
$~G~$ is not a group ). We have that $~G~$ is a monoid.\hfill$\b$

{\bf Theorem 2.1.}~({\bf theorem of classification })~{\it Let
$~(G,\alpha, \beta, \mu, \iota ; G_{0})~$ be a groupoid of type
$~(4;2).~$ Then the groupoid $~G~$ is isomorphic with one from the
following groupoids:

$(i)~~~ G~\cong~\{~e~\}\coprod {\bf Z}_{3};~~~(ii)~~~~G~\cong~{\bf
Z}_{2} \coprod {\bf Z}_{2};~~~ (iii)~~~ G~\cong~{\cal
F}_{(4;2)}({\bf R}^{2})~$

where $~{\cal F}_{(4;2)}({\bf R}^{2})~$ is the groupoid of saltus
functions defined on the axes of coordinates in a plane.}

{\bf Proof.}~Let the groupoid $~(G,\alpha,\beta,\mu,\iota;G_{0})~$
of type $~(4;2),~$  where $~G=\{~a_{1}, a_{2}, a_{3}, a_{4}~\}~$~
and $~ G_{0}=~\{~a_{1}, a_{2}~\}.~$

For the surjections $~\alpha, \beta : G~\to~G_{0}~$ with properties $~\alpha (a_{k})=\beta(a_{k})=a_{k},~$ for $~k=1,2~$ exists the following situations given by the tables :\\
$$\begin{array}{cccc} {\bf Case 1.}&
\begin{array}{|r|c|c|c|c|}\hline
a & a_{1} & a_{2} & a_{3} & a_{4}\\ \hline \alpha(a) & a_{1} &
a_{2} & a_{2} & a_{2}\cr \hline \beta(a)  & a_{1} & a_{2} & a_{2}
& a_{2}\cr \hline
\end{array}&
\ \ {\bf Case 2.}&
\begin{array}{|r|c|c|c|c|}\hline
a & a_{1} & a_{2} & a_{3} & a_{4}\\ \hline \alpha(a) & a_{1} &
a_{2} & a_{1} & a_{2}\cr \hline \beta(a)  & a_{1} & a_{2} & a_{1}
& a_{2}\cr \hline
\end{array}\end{array}$$

$${\bf Case 3.}\ \ \ \begin{array}{|r|c|c|c|c|}\hline
a & a_{1} & a_{2} & a_{3} & a_{4}\\ \hline \alpha(a) & a_{1} &
a_{2} & a_{1} & a_{2}\cr \hline \beta(a)  & a_{1} & a_{2} & a_{2}
& a_{1}\cr \hline
\end{array}$$

In the {\bf Case 1}, for the injective map $~\iota : G~\to~G~$
with property $~\iota(a_{k}) = a_{k},~$ for $~k=1,2~$ exists the
following cases given by the tables:
$$\begin{array}{cccc}
{\bf Case 1.1.}&
\begin{array}{|r|c|c|c|c|}\hline
a & a_{1} & a_{2} & a_{3} & a_{4}\\ \hline \iota(a) & a_{1} &
a_{2} & a_{3} & a_{4}\cr \hline
\end{array}&
\ \ \ {\bf Case 1.2.} &\begin{array}{|r|c|c|c|c|}\hline a & a_{1}
& a_{2} & a_{3} & a_{4}\\ \hline \iota(a) & a_{1} & a_{2} & a_{4}
& a_{3}\cr \hline
\end{array}\end{array}$$

In the {\bf Case 1} the set $~G_{(2)}~$ of composable pairs of the
groupoid $~(G;G_{0})~$ is $G_{(2)}\!=\!\{(a_{1},a_{1});
(a_{2},a_{2});(a_{2},a_{3});(a_{2},a_{4});(a_{3},a_{2});(a_{3},a_{3});(a_{3},a_{4});
(a_{4},a_{2});\linebreak(a_{4},a_{3});(a_{4},a_{4})\}.$

{\bf Case 1.1.} Using the fact that $~\mu(a_{k},a_{k})=a_{k}~$ for
$~k=1,2~$ and the properties of the functions $~\alpha, \beta~$
and $~\iota~,$ the structure table of $~\mu~$ in the situation of
{\bf Case 1.1} is the following :

$$\begin{array}{|r|c|c|c|c|}\hline
\mu   & a_{1} & a_{2} & a_{3}            & a_{4}\\ \hline a_{1} &
a_{1} &       &                  &                \cr \hline a_{2}
&       & a_{2} & a_{3}            & a_{4}           \cr \hline
a_{3} &       & a_{3} & a_{4}            & \mu(a_{3},a_{4}) \cr
\hline a_{4} &       & a_{4} & \mu(a_{4},a_{3}) & a_{3}
\cr \hline
\end{array}$$

In this case, it follows that for $~\mu~$ we have the following
situations: {\bf Cases 1.1.1 - 1.1.9} obtained by taking
$~\mu(a_{3},a_{4})\in \{~a_{2}, a_{3}, a_{4}~\}~$ and
$~\mu(a_{4},a_{3})\in \{~a_{2}, a_{3}, a_{4}~\}.~$

If we introduce the initial data in each situation of {\bf Case
1.1.1-1.1.9} and apply the program $~BGroidAP1,~$ we obtain that
$~(G;G_{0})~$ is not a groupoid.

{\bf Case 1.2.} Using the fact that $~\mu(a_{k},a_{k})=a_{k}~$ for
$~k=1,2~$ and the properties of the functions $~\alpha, \beta~$
and $~\iota~,$ the structure table of $~\mu~$ in the situation of
{\bf Case 1.2} is the following :

$$\begin{array}{|r|c|c|c|c|}\hline
\mu   & a_{1} & a_{2} & a_{3}            & a_{4}\\ \hline a_{1} &
a_{1} &       &                  &                \cr \hline a_{2}
&       & a_{2} & a_{3}            & a_{4}           \cr \hline
a_{3} &       & a_{3} & \mu(a_{3},a_{3}) & a_{2}           \cr
\hline a_{4} &       & a_{4} & a_{2}            & \mu(a_{4},a_{4})
\cr \hline
\end{array}$$

In this case, it follows that for $~\mu~$ we have the following
situations: {\bf Cases 1.2.1 - 1.2.9} obtained by taking
$~\mu(a_{3},a_{3})\in \{~a_{2}, a_{3}, a_{4}~\}~$ and
$~\mu(a_{4},a_{4})\in \{~a_{2}, a_{3}, a_{4}~\}.~$

If we introduce the initial data in each situation of {\bf Case
1.1.1-1.1.9} and apply the program $~BGroidAP1,~$ we obtain that
$~(G;G_{0})~$ is a groupoid when $~\mu(a_{3},a_{3})=a_{4}~$ and
$~\mu(a_{4},a_{4})=a_{3} ~$  ( in the other cases, $~G~$ is not a
groupoid ). For this situation, the structure functions $~\alpha,
\beta, \iota ~$ and $~\mu~$ are given by the tables:

$$\begin{array}{|r|c|c|c|c|}\hline
a & a_{1} & a_{2} & a_{3} & a_{4}\\ \hline \alpha(a) & a_{1} &
a_{2} & a_{2} & a_{2}\cr \hline \beta(a)  & a_{1} & a_{2} & a_{2}
& a_{2}\cr \hline \iota(a) & a_{1} & a_{2} & a_{4} & a_{3}\cr
\hline
\end{array}~~~~~\begin{array}{|r|c|c|c|c|}\hline
\mu   & a_{1} & a_{2} & a_{3}            & a_{4}\\ \hline a_{1} &
a_{1} &       &                  &                \cr \hline a_{2}
&       & a_{2} & a_{3}            & a_{4}           \cr \hline
a_{3} &       & a_{3} & a_{4}            & a_{2}           \cr
\hline a_{4} &       & a_{4} & a_{2}            & a_{3} \cr \hline
\end{array}$$

It is easy to prove that this groupoid is the union of two groups
and it is isomorphic with $~\{~e~\}\coprod {\bf Z}_{3}.~$ Hence,
$~G\cong \{~e~\}\coprod {\bf Z}_{3}.~$

{\bf Case 2.} Applying the same methode as in the {\bf Case 1}, we
obtain that the groupoid $~G~$ is the disjoint union $~G_{1}\cup
G_{2},~$ where $~G_{1}=\{~a_{1}, a_{3}~\}~$ such that
$~a_{3}^{2}=a_{1}~$ and $~G_{2}=\{~a_{2}, a_{4}~\}~$ such that
$~a_{4}^{2}=a_{2}.~$

For this groupoid, the structure functions $~\alpha, \beta, \iota
~$ and $~\mu~$ are given by the tables:
$$\begin{array}{|r|c|c|c|c|}\hline
a         & a_{1} & a_{2} & a_{3} & a_{4}\\ \hline \alpha(a) &
a_{1} & a_{2} & a_{1} & a_{2}\cr \hline \beta(a)  & a_{1} & a_{2}
& a_{1} & a_{2}\cr \hline \iota(a)  & a_{1} & a_{2} & a_{3} &
a_{4}\cr \hline
\end{array}~~~~~\begin{array}{|r|c|c|c|c|}\hline
\mu   & a_{1} & a_{2} & a_{3}            & a_{4}\\ \hline a_{1} &
a_{1} &       &  a_{3}           &                \cr \hline a_{2}
&       & a_{2} &                  & a_{4}           \cr \hline
a_{3} & a_{3} &       & a_{1}            &            \cr \hline
a_{4} &       & a_{4} &                  & a_{2} \cr \hline
\end{array}$$

It is easy to prove that this groupoid isomorphic with $~{\bf
Z}_{2}\coprod {\bf Z}_{2}.~$ Hence, $~G\cong {\bf Z}_{2}\coprod
{\bf Z}_{2}.~$

{\bf Case 3.} Similarly, we obtain that the structure functions of
the groupoid $~G~$ are given by the tables:
$$\begin{array}{|r|c|c|c|c|}\hline
a         & a_{1} & a_{2} & a_{3} & a_{4}\\ \hline \alpha(a) &
a_{1} & a_{2} & a_{1} & a_{2}\cr \hline \beta(a)  & a_{1} & a_{2}
& a_{2} & a_{1}\cr \hline \iota(a)  & a_{1} & a_{2} & a_{4} &
a_{3}\cr \hline
\end{array}~~~~~\begin{array}{|r|c|c|c|c|}\hline
\mu   & a_{1} & a_{2} & a_{3}            & a_{4}\\ \hline a_{1} &
a_{1} &       &  a_{3}           &                \cr \hline a_{2}
&       & a_{2} &                  & a_{4}           \cr \hline
a_{3} &       & a_{3} &                  & a_{1}           \cr
\hline a_{4} & a_{4} &       & a_{2}            &
\cr \hline
\end{array}$$

It is easy to prove that this groupoid isomorphic with $~{\cal
F}_{(4;2)}({\bf R}^{2}),~$  see Example . Hence, $~G\cong {\cal
F}_{(4;2)}({\bf R}^{2}).~$
 We observe that in this case , $~G~$ is not a group bundle.\hfill$\b$

For more details concerning the program $~BGroidAP1,~$ the reader
can be inform at
e-mail adress: ivan@hilbert.math.uvt.ro.\\[0.2cm]

\begin{center}
{\bf References}
\end{center}

{\bf [1].} R. Brown, {\it From Groups to Groupoids : a brief
survey}. Bull. London Math. Soc.,{\bf 19}, 1987, 113-134.

{\bf [2].} A. Coste, P. Dazord and A. Weinstein, {\it Groupoides
symplectiques}. Publ. Dept. Math. Lyon, 2/A, 1987,1-62.

{\bf [3].} P. J. Higgins, {\it Notes on categories and groupoids}.
Von Nostrand Reinhold, London, 1971.

{\bf [4].} Gh. Ivan, {\it Cayley Theorem for Monoidoids}. Glasnick
Matematicki, Vol. 31(51), 1996, 73-82.

{\bf [5].} Gh. Ivan, {\it Algebraic constructions of Brandt
groupoids}. Proceedings of the Algebra Symposium, " Babe\c s-
Bolyai" University, Cluj, 2002, 69-90.

{\bf [6].} M. V. Karasev, {\it Analogues of objects of Lie groups
theory for nonlinear Poisson brackets}. Math. U.S.S.R. Izv., 26,
1987, 497-527.

{\bf [7].} K. Mackenzie, {\it Lie groupoids and Lie algebroids in
differential geometry}. London Math. Soc., Lectures Notes Series,
{\bf 124}, Cambridge Univ.Press, 1987.

{\bf [8].} A. Weinstein, {\it  Groupoids: Unifying Internal and External Symmetries}. Notices Amer. Math. Soc., {\bf 43}, 1996, 744 - 752.\\[0.4cm]

Author's address:\\

West University of Timi\c soara\\
\hspace*{0.7cm} Department of Mathematics \\
\hspace*{0.7cm} 4, Bd. V. P{\^a}rvan, 1900, Timi\c soara , Romania\\
\hspace*{0.7cm} E-mail : ivan@hilbert.math.uvt.ro\\

West University of Timi\c soara\\
\hspace*{0.7cm} Department of Mathematics. Seminarul de  Algebr\u a\\
\hspace*{0.7cm}4, Bd. V. P{\^a}rvan, 1900, Timi\c soara, Romania

\end{document}